\documentclass[10pt]{amsart}

\usepackage{amsmath,amsfonts,amssymb,amsthm,amscd,epsfig}

\def\C{\mathbb{C}}
\def\Z{\mathbb{Z}}

\def\R{{\mathbb{R}}}
\def\id{\mathbf{1}}

\def\g{\ensuremath{\mathfrak{g}}}

\def\v{\mathbf{v}}
\def\w{\mathbf{w}}
\def\e{\mathbf{e}}
\def\U{\mathbf{U}}
\def\B{\mathbf{B}}
\def\F{\mathcal{F}}
\def\M{\mathfrak{M}}
\def\L{\mathfrak{L}}
\def\B{\mathfrak{B}}
\def\d{{\mathbf{d}}}
\def\t{\textrm{top}}
\def\a{\mathbf{a}}
\def\ke{{\tilde e}}
\def\kf{{\tilde f}}
\def\P{\mathcal{P}}
\def\pr{\mathrm{pr}}

\newcommand{\comment}[1]{}

\DeclareMathOperator{\Hom}{Hom}

\DeclareMathOperator{\End}{End}

\DeclareMathOperator{\inc}{in}
\DeclareMathOperator{\out}{out}
\DeclareMathOperator{\ord}{ord}

\DeclareMathOperator{\wt}{wt}

\DeclareMathOperator{\Spec}{Spec}
\DeclareMathOperator{\Id}{Id}

\newtheorem{theo}{Theorem}[section]
\newtheorem{prop}[theo]{Proposition}
\newtheorem{lem}[theo]{Lemma}
\newtheorem{cor}[theo]{Corollary}

\newtheorem*{rem*}{Remark}
\newtheorem{rem}[theo]{Remark}
\numberwithin{equation}{section}

\begin{document}
\title{On two geometric constructions of $U(\mathfrak{sl}_n)$ and its representations}
\author{Alistair Savage}
\address{University of Ottawa, Ottawa, Ontario, Canada}
\thanks{This research was supported in part by the Natural
Sciences and Engineering Research Council (NSERC) of Canada and
was partially conducted while at the Max-Planck-Institut f\" ur
Mathematik} \subjclass[2000]{16G20,17B45}
\date{May 3, 2005}

\begin{abstract}
Ginzburg and Nakajima have given two different geometric
constructions of quotients of the universal enveloping algebra of
$\mathfrak{sl}_n$ and its irreducible finite-dimensional highest
weight representations using the convolution product in the
Borel-Moore homology of flag varieties and quiver varieties
respectively. The purpose of this paper is to explain the precise
relationship between the two constructions. In particular, we show
that while the two yield different quotients of the universal
enveloping algebra, they produce the same representations and the
natural bases which arise in both constructions are the same.  We
also examine how this relationship can be used to translate the
crystal structure on irreducible components of quiver varieties,
defined by Kashiwara and Saito, to a crystal structure on the
varieties appearing in Ginzburg's construction, thus recovering
results of Malkin.
\end{abstract}

\maketitle

\section*{Introduction}
The universal enveloping algebra of $\mathfrak{sl}_n$ and its
finite-dimensional highest weight representations have been
constructed geometrically in two different ways by Ginzburg
\cite{G91} and Nakajima \cite{N98} (Nakajima's construction works
for more general Kac-Moody algebras). Both constructions use a
convolution product in homology.  In Ginzburg's construction, the
varieties involved are flag varieties and their cotangent bundles
while in Nakajima's construction they are varieties attached to the
quiver (oriented graph) whose underlying graph is the Dynkin graph
of $\mathfrak{sl}_n$.  Both realizations produce a natural basis of
the representations given by the fundamental classes of the
irreducible components of the varieties involved. In \cite{N94}
Nakajima conjectured a specific relationship between the two
varieties and this conjecture was later proved by Maffei \cite{M00}.
In the current paper we review this relationship and use it to
examine the representation theoretic constructions in the two
settings and show that while the quotients of the universal
enveloping algebra obtained are different, there is a natural
homomorphism between the two and the natural bases in
representations produced by the two constructions are in fact the
same.  Nakajima's construction using the convolution product was in
fact motivated by Ginzburg's construction and thus it is not
surprising that we find that the quiver variety construction is in
some sense a generalization of the flag variety construction to
arbitrary (simply-laced) type.  It was certainly expected by experts
that the two bases obtained are the same.  However, the author is
not aware of a proof in the literature of the coincidence of the two
bases and the precise relationship between the different
constructions of the universal enveloping algebra (which are, in
fact, slightly different in the two cases).

Finally, we use the relation between the two constructions to define
the structure of a crystal graph on the irreducible components of
the Spaltenstein varieties appearing in Ginzburg's construction by
analogy with the already existing theory for quiver varieties
developed by Kashiwara and Saito.  In doing this, we recover the
crystal structure on irreducible components of Spaltenstein
varieties introduced by Malkin in \cite{Mal02}. We now explain the
contents of the paper in some detail.

Fix a positive integer $d$ and let
\[
\F = \{0 = F_0 \subset F_1 \subset \dots \subset F_n = \C^d\}
\]
be the set of all $n$-step flags in $\C^d$.  Let $N = \{x \in
\End(\C^d)\ |\ x^n = 0 \}$.  The cotangent bundle to $\F$ is
isomorphic to
\[
M = \{(x,F) \in N \times \F\ |\ x(F_i) \subset F_{i-1}\}.
\]
We have the natural projection $\mu : M \to N$ and for $x \in N$
we define
\begin{gather*}
Z = M \times_N M = \{(m_1,m_2) \in M \times M\ |\ \mu(m_1) =
\mu(m_2)\}, \\
\F_x = \mu^{-1}(x).
\end{gather*}
Using the convolution product (see Section~\ref{sec:convolution}),
we give the top-dimensional Borel-Moore homology $H_\t(Z)$ the
structure of an algebra and $H_\t(\F_x)$ the structure of a module
over this algebra.  Let $I_d$ be the annihilator of
$(\C^n)^{\otimes d}$, a two-sided ideal of finite codimension in
the enveloping algebra $U(\mathfrak{sl}_n)$.  Here $\C^n$ is the
natural $\mathfrak{sl}_n$-module. Then in \cite{CG,G91} it is
shown that $H_\t(Z) \cong U(\mathfrak{sl}_n)/I_d$ and that under
this isomorphism, $H_\t(\F_x)$ is the irreducible highest weight
$\mathfrak{sl}_n$-module of highest weight $w_1 \omega_1 + \dots +
w_{n-1} \omega_{n-1}$ where the $\omega_i$ are the fundamental
weights of $\mathfrak{sl}_n$ and $w_i$ is the number of $(i \times
i)$-Jordan blocks in the Jordan normal form of $x$.

Now, in \cite{N98}, Nakajima constructs the same representations
in a similar way using a convolution product in the homology of
quiver varieties.  In \cite{M00}, Maffei showed that the varieties
of Nakajima's construction are isomorphic to the following.  Let
$S_x$ be a transversal slice in $N$ to the $GL(\C^d)$-orbit
through $x$ (see Section~\ref{sec:isom}).  Then let
\begin{gather*}
M' = \mu^{-1}(S_x), \\
Z' = M' \times_{S_x} M'.
\end{gather*}
Then, translated via the isomorphism of \cite{M00}, a result of
\cite{N98} is that, under the convolution product we have
$H_\t(Z') \cong U(\mathfrak{sl}_n)/J$ and $H_\t(M_x)$ is the same
irreducible highest weight module as in Ginzburg's construction
(see Theorems~\ref{thm:nak-hom} and~\ref{thm:nak-quotient}).  Here
$J$ is a certain ideal of finite codimension in
$U(\mathfrak{sl}_n)$ that is different from $I_d$ in general. Thus
the two constructions yield different quotients of the universal
enveloping algebra but the same representation.

Since $Z' \subset Z$ and $M' \subset M$, we have a natural
restriction with support morphism $H_\t(Z) \to H_\t(Z')$.  The
main result of this paper (see Theorem~\ref{thm:quiver-flag}) is
that the following diagram is commutative
\[
\begin{CD}
H_\t(Z) \otimes H_\t(\F_x) @>>> H_\t(Z') \otimes H_\t(\F_x)
@>\cong>> \oplus_{\v^1,\v^2} H_\t(Z(\v^1,\v_2;\w)) \otimes
\oplus_\v H_\t(\L(\v,\w))
\\
@VVV  @VVV  @VVV \\
H_\t(\F_x) @>>> H_\t(\F_x) @>\cong>> \oplus_\v H_\t(\L(\v,\w))
\end{CD}
\]
Here the rightmost term in each row involves the Nakajima quiver
varieties (see Section~\ref{sec:nak} for definitions).  We are
also able to conclude that the natural bases of representations
produced by both Ginzburg's and Nakajima's constructions coincide.
We thus obtain a precise relation between the two approaches.

Recently, a relation has been established between a construction
closely related to that of Ginzburg and another geometric approach
of Mirkovi\' c-Vilonen in terms of the affine Grassmannian
\cite{BGV04}. It would be interesting to examine the connection
between the quiver variety and Mirkovi\' c-Vilonen realizations of
finite-dimensional representations of Lie algebras.

The organization of the paper is as follows.  In
Sections~\ref{sec:prelim} and \ref{sec:convolution} we recall the
definition of $\mathfrak{sl}_n$ and the convolution product in
Borel-Moore homology.  In Sections~\ref{sec:ginz-def} and
\ref{sec:nak} we review Ginzburg's and Nakajima's constructions of
$U(\mathfrak{sl}_n)$ and its representations.  Then in
Section~\ref{sec:isom} we describe the precise relationship
between the two constructions.  Finally, in
Section~\ref{sec:crystal} we define the structure of a crystal on
the irreducible components of $\F_x$.

The author would like to thank O. Schiffmann for many useful
discussions and suggestions and K. McGerty for very helpful
comments on the properties of the convolution product.

%%%%%%%%%%%%%%%%%%%%%%%%%%%%%%%%%%%%%%%%%%%%%%%%%%%%%%%%%%%%%%%%%%%%%%%

\section{Preliminaries}
\label{sec:prelim}

Let $\g=\mathfrak{sl}_n$ be the Lie algebra of type $A_{n-1}$.
Then \g\ is the space of all traceless $n \times n$ matrices.
Let $\{e_k,f_k\}_{k=1}^{n-1}$ be the set of Chevalley generators.
The Cartan subalgebra $\mathfrak{h}$ is spanned by the matrices
\[
h_k = e_{k,k} - e_{k+1,k+1},\ 1 \le k \le n-1,
\]
where $e_{k,l}$ is the matrix with a one in entry $(k,l)$ and
zeroes everywhere else.  Thus the dual space $\mathfrak{h}^*$ is
spanned by the simple roots
\[
\alpha_k = \epsilon_k - \epsilon_{k+1},\ 1\le k \le n-1,
\]
where $\epsilon_k(e_{l,l}) = \delta_{kl}$ and the fundamental
weights are given by
\[
\omega_k = \epsilon_1 + \dots + \epsilon_k,\ 1 \le k \le n-1.
\]
Consider a dominant weight $\w = w_1 \omega_1 + \dots + w_{n-1}
\omega_{n-1}$.  Then
\[
\w = \lambda_1 \epsilon_1 + \dots + \lambda_{n-1} \epsilon_{n-1}
\]
where $\lambda_k = w_k + \dots + w_{n-1}$ and so $\w$ corresponds to
a partition $\lambda(\w) = (\lambda_1 \ge \dots \ge \lambda_{n-1})$.
We say that a highest weight $\w$ is a partition of $d$ if
$|\lambda(\w)| = \lambda_1 + \dots + \lambda_{n-1} = d$ or,
equivalently, if $\sum_{k=1}^n k w_k = d$.

%%%%%%%%%%%%%%%%%%%%%%%%%%%%%%%%%%%%%%%%%%%%%%%%%%%%%%%%%%%%%%%%%%%%%%%

\section{Convolution algebra in homology}
\label{sec:convolution}

In this section we give a brief overview of the convolution algebra
in homology. The reader interested in further details should consult
\cite{CG}.

In this paper $H_*(Z)$ will denote the Borel-Moore homology with
$\C$-coefficients of a locally-compact space $Z$.  Thus, by
definition, if $Z$ is a closed subset of a smooth, oriented
manifold $M$, then
\[
H_k(Z) = H^{\dim_\R M - k}(M,M\backslash Z).
\]
If $Z$ and $Z'$ are closed subsets of a smooth variety $M$, we
have a $\cup$-product map
\[
H^k(M, M\backslash Z) \times H^l(M, M\backslash Z') \to H^{k+l}(M,
M \backslash (Z \cap Z')).
\]
Thus we construct the intersection pairing in Borel-Moore homology
\[
\cap : H_k(Z) \times H_l(Z') \to H_{k+l-d}(Z \cap Z'),\quad d=
\dim_\R M.
\]

Let $M_1$, $M_2$ and $M_3$ be smooth, oriented manifolds and
$p_{kl} : M_1 \times M_2 \times M_3 \to M_k \times M_l$ be the
obvious projections.  Let $Z \subset M_1 \times M_2$ and $Z'
\subset M_2 \times M_3$ be closed subvarieties and assume that the
map
\[
p_{13} : p_{12}^{-1}(Z) \cap p_{23}^{-1}(Z') \to M_1 \times M_3
\]
is proper and denote its image by $Z \circ Z'$.  The operation of
convolution
\[
\star: H_k(Z) \times H_l(Z') \to H_{k+l-d}(Z \circ Z'),\quad
d=\dim_\R M_2,
\]
is defined by
\[
c \star c' = (p_{13})_* (p_{12}^*c \cap p_{23}^*c'),
\]
where $p_{12}^*c$ means $c \boxtimes [M_3]$, etc.

Now, let $M$ be a smooth manifold and $\mu: M \to N$ be a proper
morphism.  Let
\[
Z = M \times_N M = \{(m_1,m_2) \in M \times M\ |\ \mu(m_1) =
\mu(m_2)\} \subset M \times M.
\]
Then $Z \circ Z = Z$ and so convolution makes $H_*(Z)$ a
finite-dimensional associative $\C$-algebra with unit.

For $x \in N$, let $M_x = \mu^{-1}(x)$.  We also identify $M_x$
with the variety $M_x \times \text{pt}$.  Then setting $M_1 = M_2
= M$ and $M_3 = \text{pt}$, we have $Z \circ M_x = M_x$ and
convolution makes $H_*(M_x)$ a $H_*(Z)$-module.

%%%%%%%%%%%%%%%%%%%%%%%%%%%%%%%%%%%%%%%%%%%%%%%%%%%%%%%%%%%%%%%%%%%%%%%%%%%%

\section{Ginzburg's construction}
\label{sec:ginz-def}

We recall here Ginzburg's construction of the enveloping algebra
$U(\mathfrak{sl}_n)$ and its irreducible highest weight
representations.  Proofs omitted here can be found in \cite{G91} or
\cite{CG}.

Fix an integer $d \ge 1$.  Let
\[
\F = \{ 0 = F_0 \subset F_1 \subset \dots \subset F_n = \C^d\}
\]
be the set of all $n$-step partial flags in $\C^d$.  The space $\F$
is a disjoint union of smooth compact manifolds with connected
components parameterized by compositions
\[
\d = (d_1 + d_2 + \dots + d_n = d), \quad d_i \in \Z_{\ge 0}.
\]
The connected component of $\F$ corresponding to $\d$ is
\[
\F_\d = \{F = (0 = F_0 \subset \dots \subset F_n = \C^d)\ |\ \dim
F_i/F_{i-1} = d_i\},
\]
and
\[
\dim_\C \F_\d = \frac{d!}{d_1! d_2! \cdots d_n!}.
\]

Let
\[
N = \{x \in \End(\C^d)\ |\ x^n=0\}.
\]
Then
\[
T^*\F \cong M = \{(x,F) \in N \times \F\ |\ x(F_i) \subset
F_{i-1},\ 1 \le i \le n\}.
\]
The above decomposition of $\F$ yields a decomposition of $M$ given
by $M = \sqcup_\d M_\d$ where $M_\d = T^*\F_\d$ for an $n$-step
composition $\d$ of $d$.

The natural projections give rise to the diagram
\[
N \stackrel{\mu}{\longleftarrow} M
\stackrel{\pi}{\longrightarrow} \F.
\]
We have a natural action of $GL_d(\C)$ on $\F$, $N$ (by conjugation)
and $M$ and the projections commute with this action.

For $x \in N$, let $\F_x = \mu^{-1}(x)$.  It has connected
components $\F_{\d,x}$ given by $\F_{\d,x} = \F_\d \cap \F_x$.
Define
\[
Z = M \times_N M = \{(m_1,m_2) \in M \times M\ |\ \mu(m_1) =
\mu(m_2) \} \subset M \times M.
\]
We use the convention that under the isomorphism
\[
T^*\F \times T^*\F \cong T^*(\F \times \F),
\]
the standard symplectic form on the right hand side corresponds to
$\omega_1 - \omega_2$ where $\omega_1$ and $\omega_2$ are the
symplectic forms on the first and second factors of the left hand
side respectively.

\begin{prop}
The variety $Z$ is the union of the conormal bundles to the
$GL_d(\C)$-orbits in $\F \times \F$.  The closures of these
conormal bundles are precisely the irreducible components of $Z$.
\end{prop}

\begin{prop}
We have $Z \circ Z = Z$.  Thus $H_*(Z)$ is an associative algebra
with unit and $H_*(\F_x)$ is an $H_*(Z)$-module for any $x \in
N$.
\end{prop}

\begin{prop}
All irreducible components of $Z$ contained in $M_{\d^1} \times
M_{\d^2}$ are half dimensional.  That is, they have complex
dimension
\begin{gather*}
\frac{1}{2} \dim_\C(M_{\d^1} \times M_{\d^2}) = \frac{1}{2} \left(
2\frac{d^1!}{d^1_1! d^1_2! \cdots d^1_n!} + 2\frac{d^2!}{d^2_1!
d^2_2! \cdots d^2_n!} \right) \\
= \frac{d^1!}{d^1_1! d^1_2! \cdots d^1_n!} + \frac{d^2!}{d^2_1!
d^2_2! \cdots d^2_n!}.
\end{gather*}
\end{prop}

Let $H_\t(Z)$ be the vector subspace of $H_*(Z)$ spanned by the
fundamental classes of the irreducible components of $Z$ and let
$H_\t(\F_x)$ be the vector subspace of $H_*(\F_x)$ spanned by the
fundamental classes of the irreducible components of $\F_x$.

\begin{prop}
The homology group $H_\t(Z)$ is a subalgebra of $H_*(Z)$ and
$H_\t(\F_x)$ is an $H_\t(Z)$-stable subspace of $H_*(\F_x)$.
\end{prop}

Now, for a composition $\d$ we have the diagonal subvariety $\Delta
\subset \F_\d \times \F_\d$ which is a $GL_d(\C)$-orbit. We define
\[
H_k = \sum_\d (d_k - d_{k+1}) [T^*_\Delta (\F_\d \times \F_\d)],
\]
where $T_O^*(\F_\d \times \F_\d)$ denotes the conormal bundle
to a $GL_d(\C)$-orbit $O \subset \F_\d \times \F_\d$.
Note that under the sign convention for the symplectic form
mentioned above, the conormal bundle $T^*_\Delta (\F_\d \times
\F_\d)$ is the diagonal in $T^*\F_\d \times T^*\F_\d$.

Now, for a composition $\d = (d_1 + \dots + d_n)$ and $1 \le k \le
n-1$, let
\begin{align*}
\d_k^+ &= d_1 + \dots + d_{k-1} + (d_k+1) + (d_{k+1}-1) + d_{k+2}
+ \dots + d_n, \\
\d_k^- &= d_1 + \dots + d_{k-1} + (d_k-1) + (d_{k+1}+1) + d_{k+2}
+ \dots + d_n,
\end{align*}
provided that these are compositions (that is, all terms are $\ge
0$). Otherwise, we define $\d_k^{\pm} = \nabla$, the ghost
composition.

If $1 \le k \le n-1$ and $\d = (d_1 + \dots + d_n)$ is a composition
such that $\d_k^+ \ne \nabla$, resp. $\d_k^- \ne \nabla$, we define
\begin{align*}
Y_{\d_k^+,\d} &= \{ (F',F) \in \F_{\d_k^+} \times \F_\d \ |\ F_l =
F'_l \ \forall \ l\ne k,\, F_k \subset F_k', \dim (F_k'/F_k) = 1 \},\\
Y_{\d_k^-,\d} &= \{ (F',F) \in \F_{\d_k^-} \times \F_\d \ |\ F_l =
F'_l \ \forall \ l\ne k,\, F_k' \subset F_k, \dim (F_k/F_k') = 1
\}.
\end{align*}
Note that each $Y_{\d_k^{\pm},\d}$ is a $GL_d(\C)$-orbit in
$\F_{\d_k^{\pm}} \times \F_\d$ of minimal dimension and thus is a
smooth closed subvariety.  Let
\begin{align}
E_k &= \sum_\d [T^*_{Y_{\d_k^+,\d}}(\F_{\d_k^+} \times \F_\d)],\\
\label{def:ginz-Fk}
F_k &= \sum_\d
(-1)^{s_k(\d_k^+,\d)}[T^*_{Y_{\d_k^-,\d}}(\F_{\d_k^-} \times
\F_\d)],
\end{align}
where $s_k(\d_k^+,\d) = \frac{1}{2}\left( \dim_\C M_{\d_k^+} -
\dim_\C M_\d \right)$.

\begin{theo}[\cite{G91}]
\label{thm:ginz}
The map
\[
e_k \mapsto E_k,\ f_k \mapsto F_k,\ h_k \mapsto H_k,
\]
extends to a surjective algebra homomorphism $U(\mathfrak{sl}_n)
\twoheadrightarrow H_\t(Z)$.  Under this homomorphism,
$H_\t(\F_x)$ is the irreducible highest weight module of highest
weight $w_1 \omega_1 + \dots + w_{n-1} \omega_{n-1}$ where
$\omega_i$ are the fundamental weights and $w_i$ is the number of
($i \times i$)-Jordan blocks in the Jordan normal form of $x$.
\end{theo}

\begin{rem}
Note that the sign appearing in \eqref{def:ginz-Fk} does not
appear in \cite{CG,G91}.  This arises from the fact that
Theorem~2.7.26 (iii) in \cite{CG} should read $[Z_{12}] \star
[Z_{23}] = (-1)^{\dim F} \chi(F) \cdot [Z_{13}]$ (see \cite[Lemma
8.5]{N98}).
\end{rem}

Let $I_d$ be the annihilator of $(\C^n)^{\otimes d}$, a two-sided
ideal of finite codimension in the enveloping algebra
$U(\mathfrak{sl}_n)$.  Here $\C^n$ is the natural
$\mathfrak{sl}_n$-module.

\begin{theo}[{\cite[Proposition~4.2.5]{CG}}]
\label{thm:ginz-quotient}
The homomorphism of
Theorem~\ref{thm:ginz} yields an algebra isomorphism
\[
U(\mathfrak{sl}_n)/I_d \cong H_\t(Z).
\]
\end{theo}

It is known that the simple $\mathfrak{sl}_n$ modules that occur
with non-zero multiplicity in the decomposition of
$(\C^n)^{\otimes d}$ are precisely those modules whose highest
weight is a partition of $d$.

%%%%%%%%%%%%%%%%%%%%%%%%%%%%%%%%%%%%%%%%%%%%%%%%%%%%%%%%%%%%%%%%%%%%%%

\section{Nakajima's construction}
\label{sec:nak}

In this section, we will review the description of the quiver
varieties presented in \cite{N98}.  Further details may be found in
\cite{N94} and \cite{N98}.  We only discuss the case corresponding
to the Lie algebra $\mathfrak{sl}_n$.  Note that we use a different
stability condition that the one used in \cite{N94} and \cite{N98}
and so our definitions differ slightly from the ones that appear
there. One can translate between the two stability conditions by
taking transposes of the maps appearing in the definitions of the
quiver varieties. See \cite{N96} for a discussion of various choices
of stability condition.

\begin{figure}
\centering \epsfig{file=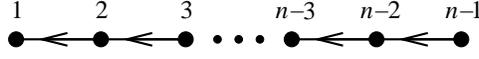,width=0.5\textwidth}
\caption{The quiver of type $A_{n-1}$. \label{fig:quiver_an}}
\end{figure}

As before, let $\mathfrak{g}=\mathfrak{sl}_n$ be the simple Lie
algebra of type $A_{n-1}$. Let $I=\{1,\dots,n-1\}$ be the set of
vertices of the Dynkin graph of $\g$ with the set of oriented
edges given by
\begin{gather*}
H=\{h_{k,l} \ |\ k,l \in I,\ |k-l|=1\}.
\end{gather*}
For two adjacent vertices $k$ and $l$, $h_{k,l}$ is the oriented
edge from vertex $k$ to vertex $l$.  We denote the outgoing and
incoming vertices of $h\in H$ by $\out(h)$ and $\inc(h)$
respectively.  Thus $\out(h_{k,l}) = k$ and $\inc(h_{k,l})=l$.
Define the involution $\bar{\ } : H \to H$ as the function that
interchanges $h_{k,l}$ and $h_{l,k}$.  Fix the orientation $\Omega
= \{h_{k,k-1}\ |\ 2 \le k \le n-1 \}$.  We picture this quiver as
in Figure~\ref{fig:quiver_an}.

Let $V = \bigoplus_{k \in I} V_k$ and $W = \bigoplus_{k \in I}
W_k$ be two finite dimensional complex $I$-graded vector spaces with graded
dimensions
\begin{align*}
\v &= (\dim V_1,\dim V_2,\dots,\dim V_{n-1}),\\
\w &= (\dim W_1,\dim W_2,\dots,\dim W_{n-1}).
\end{align*}
Then we define
\[
\mathbf{M}(\v,\w) = \bigoplus_{h \in H}
\Hom(V_{\out(h)},V_{\inc(h)}) \oplus \bigoplus_{k \in I}
\Hom(W_k,V_k) \oplus \bigoplus_{k \in I} \Hom(V_k,W_k).
\]
The above three components of an element of $\mathbf{M}(\v,\w)$
will be denoted by $B=(B_h)$, $i=(i_k)$ and $j=(j_k)$.  We
associate elements in the weight lattice of $\g$ to the dimensions
vectors $\v = (v_1,\dots,v_{n-1})$ and $\w=(w_1,\dots,w_{k-1})$ as
follows.
\[
\alpha_\v = \sum_{k \in I} v_k \alpha_k,\quad \omega_\w = \sum_{k
\in I} w_k \omega_k,
\]
where $\alpha_k$ and $\omega_k$ are the simple roots and
fundamental weights respectively.

Now, let
\[
G_\v = \prod_{k \in I} GL(V_k)
\]
act on $\mathbf{M}(\v,\w)$ by
\[
g(B,i,j) = (gBg^{-1},gi,jg^{-1}),
\]
where $gBg^{-1} = (B'_h) = (g_{\inc(h)}B_h g_{\out(h)}^{-1})$, $gi
= (i'_k) = (g_k i_k)$ and $jg^{-1} = (j'_k) = (j_k g_k^{-1})$. Let
$\epsilon : H \to \{\pm 1\}$ be given by
\[
\epsilon(h) = \begin{cases}
+1 & \text{if } h \in \Omega \\
-1 & \text{if } h \in {\bar \Omega}
\end{cases}.
\]
Define a map $\mu : \mathbf{M}(\v,\w) \to \bigoplus_{k \in I} \End (V_k,
V_k)$ with $k$th component given by
\[
\mu_k(B,i,j) = \sum_{h \in H\, :\, \inc(h)=k} \epsilon(h) B_h
B_{\bar h} + i_k j_k.
\]
Let $A(\mu^{-1}(0))$ be the coordinate ring of the affine
algebraic variety $\mu^{-1}(0)$ and define
\[
\M_0(\v,\w) = \mu^{-1}(0)//G = \Spec A(\mu^{-1}(0))^G.
\]
This is the affine algebro-geometric quotient of $\mu^{-1}(0)$ by
$G$.  It is an affine algebraic variety and its geometric points
are closed $G_\v$-orbits.

We say that a collection $S=(S_k)$ of subspaces $S_k \subset V_k$
is $B$-stable if $B_h(S_{\out(h)}) \subset S_{\inc(h)}$ for all $h
\in H$.  We say that a point of $\mu^{-1}(0)$ is stable if any
$B$-stable collection of subspaces $S$ containing the image of $i$
is equal to all of $V$.  We let $\mu^{-1}(0)^s$ denote the set of
stable points.

\begin{prop}
The stabilizer in $G_\v$ of any point in $\mu^{-1}(0)^s$ is
trivial.
\end{prop}

We then define
\[
\M(\v,\w) = \mu^{-1}(0)^s/G_\v,
\]
which is diffeomorphic to an affine algebraic manifold.  We know
(see \cite[Cor 3.12]{N98}) that
\[
\dim_\C \M(\v,\w) = \v \cdot (2\w - C \v),
\]
where $C$ is the Cartan matrix of $\mathfrak{sl}_n$.

For $(B,i,j) \in \mu^{-1}(0)^s$, we denote the corresponding orbit
in $\M(\v,\w)$ by $[B,i,j]$ and if the orbit through $(B,i,j)$ is
closed, we denote the corresponding point of $\M_0(\v,\w)$ by the
same notation.

We have a map
\[
\pi : \M(\v,\w) \to \M_0(\v,\w)
\]
which sends an orbit $[B,i,j]$ to the unique closed orbit
$[B_0,i_0,j_0]$ contained in the closure of $G(B,i,j)$.  Let
$\L(\v,\w) = \pi^{-1}(0)$.

\begin{prop}
The subvariety $\L(\v,\w) \subset \M(\v,\w)$ is half-dimensional
and is homotopic to $\M(\v,\w)$.
\end{prop}
Actually, under a natural symplectic form on $\M(\v,\w)$, the
subvariety $\L(\v,\w)$ is Lagrangian.  It will be useful in the
sequel to also consider the following direct construction of
$\L(\v,\w)$.  Let
\[
\Lambda(\v,\w) = \{(B,i,j) \in \mu^{-1}(0)\ |\ j=0,\, B \text{ is
nilpotent}\}
\]
where $B$ nilpotent means that there exists an $N \ge 1$ such that
for any sequence $h_1, h_2, \dots, h_N$ in $H$ satisfying
$\inc(h_k) = \out(h_{k+1})$, the composition $B_{h_N} \cdots B_{h_2}
B_{h_1} : V_{\out(h_1)} \to V_{\inc(h_N)}$ is zero.
Furthermore, define
\[
\Lambda(\v,\w)^s = \{(B,i,j) \in \Lambda(\v,\w)\ |\ (B,i,j) \in
\mu^{-1}(0)^s\}.
\]
Then we have the following Lemma.
\begin{lem}
We have
\[
\L(\v,\w) = \Lambda(\v,\w)^s/G_\v.
\]
\end{lem}

If $V'=(V'_k)$ is a collection of subspaces of $V=(V_k)$, we have
a natural inclusion map $\M_0(\v',\w) \hookrightarrow
\M_0(\v,\w)$. Thus, for vector spaces $V^1,V^2,W$, we can consider
the projections $\pi : \M(\v^k,\w) \to \M_0(\v^k,\w)$ as maps to
$\M_0(\v^1+\v^2,\w)$.  We then define
\[
Z(\v^1,\v^2;\w) = \{(x^1,x^2) \in \M(\v^1,\w) \times \M(\v^2,\w)\
|\ \pi(x_1) = \pi(x_2)\}.
\]
Since $Z(\v^1,\v^2;\w) \circ Z(\v^2,\v^3;\w) \subset
Z(\v^1,\v^3;\w)$, we have the convolution product
\[
H_*(Z(\v^1,\v^2;\w)) \otimes H_*(Z(\v^2,\v^3;\w)) \to
H_*(Z(\v^1,\v^3;\w)).
\]
All of the irreducible components of $Z(\v^1,\v^2;\w)$ have the same dimension.
Let $H_\t(Z(\v^1,\v^2;\w))$ denote the top degree part of
$H_*(Z(\v^1,\v^2;\w))$.  It has a natural basis $\{[X]\}$ where
$X$ runs over the irreducible components of $Z(\v^1,\v^2;\w)$.

\begin{prop}
The convolution product makes the direct sum
$\bigoplus_{\v^1,\v^2} H_*(Z(\v^1,\v^2;\w))$ into an associative
algebra, and $\bigoplus_\v H_*(\L(\v,\w))$ is a left
$\bigoplus_{\v^1,\v^2} H_*(Z(\v^1,\v^2;\w))$-module.  In addition,
the top degree part $\bigoplus_{\v^1,\v^2} H_\t(Z(\v^1,\v^2;\w))$
is a subalgebra, and $\bigoplus_\v H_\t(\L(\v,\w))$ is a
$\bigoplus_{\v^1,\v^2} H_\t(Z(\v^1,\v^2;\w))$-stable submodule.
\end{prop}

Let $\Delta(\v,\w)$ denote the diagonal in $\M(\v,\w) \times
\M(\v,\w)$.  Then its fundamental class $[\Delta(\v,\w)]$ is in
$H_\t(Z(\v,\v;\w))$.  Left and right multiplication by
$[\Delta(\v,\w)]$ define projections
\begin{gather*}
[\Delta(\v,\w)] \cdot : \bigoplus_{\v^1,\v^2}
H_\t(Z(\v^1,\v^2;\w))
\to \bigoplus_{\v^2} H_\t(Z(\v,\v^2;\w)), \\
\cdot [\Delta(\v,\w)] : \bigoplus_{\v^1,\v^2}
H_\t(Z(\v^1,\v^2;\w))
\to \bigoplus_{\v^1} H_\t(Z(\v^1,\v;\w)). \\
\end{gather*}

For $k \in I$, define the \emph{Hecke correspondence}
$\B_k(\v,\w)$ to be the variety of all $(B,i,j,S)$ (modulo the
$G_\v$-action) such that $(B,i,j) \in \mu^{-1}(0)^s$ and $S$ is a
$B$-invariant subspace contained in the kernel of $j$ such that
$\dim S = \mathbf{e}^k$ where $\mathbf{e}^k$ has $k$-component
equal to one and all other components equal to zero.  We consider
$(B,i,j,S)$ as a point in $Z(\v-\mathbf{e}^k,\v;\w)$ by taking the
quotient by the subspace $S$ in the first factor.  Then
$\B_k(\v,\w)$ is an irreducible component of
$Z(\v-\mathbf{e}^k,\v;\w)$. Let $\omega : \M(\v^1,\w) \times
\M(\v^2,\w) \to \M(\v^2,\w) \times \M(\v^1,\w)$ be the map that
interchanges the two factors. Then define
\begin{gather}
E_k = \sum_{\v} [\B_k(\v,\w)] \in \bigoplus_{\v^1,\v^2}
H_\t(Z(\v^1,\v^2;\w)), \\
\label{def:nak-Fk}
F_k = \sum_{\v} (-1)^{r_k(\v,\w)} [\omega(\B_k(\v,\w))] \in
\bigoplus_{\v^1,\v^2} H_\t(Z(\v^1,\v^2;\w)), \\
H_k = \sum_{\v} \left<h_k, \omega_\w - \alpha_\v\right>
[\Delta(\v,\w)],
\end{gather}
where $r_k(\v,\w) = \frac{1}{2}(\dim \M_\C(\v - \e^k,\w) - \dim_\C
\M(\v,\w)) = -\mathbf{e}^k \cdot (\w - C\v)-1$. Here $C$ is the
Cartan matrix of $\mathfrak{sl}_n$. Note that since we are
restricting ourselves to the Lie algebra $\mathfrak{sl}_n$, the
varieties $\M(\v,\w)$ are only nonempty for a finite number of
$\v$ and so the above elements are well-defined.

\begin{theo}[\cite{N98}]
\label{thm:nak-hom} There exists a unique surjective algebra
homomorphism
\[
\Phi : U(\mathfrak{sl}_n) \twoheadrightarrow \bigoplus_{\v^1,\v^2}
H_\t(Z(\v^1,\v^2;\w))
\]
such that
\[
\Phi(h_k) = H_k, \quad \Phi(e_k) = E_k,\quad \Phi(f_k) = F_k.
\]
Under this homomorphism, $\bigoplus_\v H_\t(\L(\v,\w))$ is the
irreducible integrable highest weight module with highest weight
$\omega_\w$. The class $[\L(\mathbf{0},\w)]$ is a highest weight
vector.
\end{theo}

\begin{rem}
The result in \cite{N98} is actually in terms of the modified
universal enveloping algebra.  In the more general case of a
Kac-Moody algebra with symmetric Cartan matrix, this language is
more natural.  However, in our case of $\mathfrak{sl}_n$, since
for a fixed $\w$ the quiver varieties $\M(\v,\w)$ are non-empty
only for a finite number of $\v$, we can avoid the use of
the modified universal enveloping algebra.
\end{rem}

Let $J_\w$ be the annihilator in $U(\mathfrak{sl}_n)$ of
$\bigoplus_\v L(\omega_\w - \alpha_\v)$, where the sum is over all
$\v$ such that $\omega_\w - \alpha_\v$ is dominant integral and is a
weight of $L(\omega_\w)$.  Here $L(\lambda)$ is the irreducible
integrable highest weight representation of highest weight
$\lambda$.

\begin{theo}[{\cite[Theorem~10.15]{N98}}]
\label{thm:nak-quotient}
The homomorphism of
Theorem~\ref{thm:nak-hom} yields an algebra isomorphism
\[
U(\mathfrak{sl}_n)/J_\w \cong \bigoplus_{\v_1,\v_2}
H_\t(Z(\v^1,\v^2;\w)).
\]
\end{theo}

%%%%%%%%%%%%%%%%%%%%%%%%%%%%%%%%%%%%%%%%%%%%%%%%%%%%%%%%%%%%%%%%%%%%%%%

\section{A comparison of the two constructions}
\label{sec:isom}

We now describe the precise relationship between the constructions
of Ginzburg and Nakajima.

We begin by recalling a result of Maffei \cite{M00}.  Let $x \in
N$ and let $\{x,y,h\}$ be an $\mathfrak{sl}_2$ triple in
$GL(\C^d)$.  We define the \emph{transversal slice} to the orbit
$O_x$ of $x$ in $N$ at the point $x$ to be
\[
S_x = \{u \in N\ |\ [u-x,y]=0\}.
\]
We allow $\{0,0,0\}$ to be an $\mathfrak{sl}_2$ triple.  Thus we
have $S_0 = N$.

Now, the orbits of the action of $GL(\C^d)$ on $N$ are determined
by partitions of $d$.  Corresponding to a partition $\lambda$ is
the orbit consisting of all those matrices whose Jordan blocks
have sizes $\lambda_i$.  We let $O_\lambda$ denote the orbit
corresponding to the partition $\lambda$.

Let $\mu_\d : M_\d \to N$ denote the restriction of the map $\mu$
to $M_\d$.  Then let $\alpha = (\alpha_1 \ge \alpha_2 \ge \dots
\ge \alpha_n)$ be a permutation of $\d$ and define the partition
$\lambda_\d = 1^{\alpha_1 - \alpha_2} 2^{\alpha_2-\alpha_3} \cdots
n^{\alpha_n}$.  Then $\lambda_\d$ is a partition of $d$ and if
$(x,F) \in M_\d$, then $x \in \overline{O}_{\lambda_\d}$.
Furthermore, the map $\mu_\d : M_\d \to \overline{O}_{\lambda_\d}$
is a resolution of singularities and is an isomorphism over
$O_{\lambda_\d}$.  Define
\begin{gather*}
S_{\d,x} = S_x \cap \overline{O}_{\lambda_\d},\quad \widetilde
S_{\d,x} = \mu_\d^{-1}(S_{\d,x}) = \mu_\d^{-1}(S_x).
\end{gather*}

Now, for $\v,\w \in (\Z_{\ge 0})^{n-1}$ define $\a = \a(\v,\w) =
(a_1,\dots,a_n)$ by
\begin{gather}
\label{eq:a-def}
a_1 = w_1 + \dots + w_{n-1} - v_1,\quad a_n=v_{n-1}, \\
a_k = w_k + \dots + w_{n-1} - v_k + v_{k-1},\ 2 \le k \le n-1
\nonumber.
\end{gather}
Note that $\sum_{k=1}^n a_k = d = \sum_{k=1}^{n-1} kw_k$ and that
for a fixed $d$ and $\w$, the above map is a bijection between
$(n-1)$-tuples of integers $\v$ and $n$-tuples of integers $\a$
such that $\sum_i a_i = d$.  Furthermore, let $\M^1(\v,\w) = \pi
(\M(\v,\w))$.

\begin{theo}[\cite{M00}]
\label{thm:maffei} Let $\v$, $\w$, $d$ and $\a=\a(\v,\w)$ be as
above and let $x \in N$ be a nilpotent element of type $1^{w_1}
2^{w_2} \cdots (n-1)^{w_{n-1}}$.  Then there exists an isomorphism
$\theta : \M(\v,\w) \stackrel{\cong}{\longrightarrow} \widetilde
S_{\a,x}$ and $\theta_1 : \M^1(\v,\w)
\stackrel{\cong}{\longrightarrow} S_{\a,x}$ such that
$\theta_1(0)=x$ and the following diagram commutes.
\[
\begin{CD}
\M(\v,\w) @>\theta>> \widetilde S_{\a,x} \\
@V\pi VV                    @V \mu_\a VV  \\
\M^1(\v,\w) @>\theta_1>> S_{\a,x}
\end{CD}
\]
\end{theo}

Note that by Theorem~\ref{thm:maffei}, if we restrict $\theta$ to
$\L(\v,\w)$, we obtain an isomorphism $\L(\v,\w) \cong \F_{\a,x}$
which we will also denote by $\theta$.  This restriction is fairly
simple to describe as we now show.

We define a \emph{path} to be an ordered set of edges
$(h_1,\dots,h_N)$ such that $\inc(h_i) = \out(h_{i+1})$.  Then let
$\mathcal{P}$ be the set of all paths that head left and then
right.  That is,
\[
\mathcal{P} =
\{(h_{k,k-1},h_{k-1,k-2},\dots,h_{l+1,l},h_{l,l+1},\dots,
h_{m-1,m})\ |\ 1 \le l \le m,k \le n-1\}.
\]
For $p = (h_1,\dots,h_N) \in \mathcal{P}$, let $\inc(p) =
\inc(h_N)$ be the incoming vertex of the last edge in $p$ and let
$\out(p) = \out(h_1)$ be the outgoing vertex of the first edge in
$p$.  We define $\ord(p)$ to be the number of edges heading to the
left.  That is, $\ord(p) = \#\{h_i \in p \ |\  h_i \in \Omega\}$.
Furthermore we let $B_p = B_{h_N} \dots B_{h_1}$ be the obvious
composition of maps.

Now, for $1 \le m \le k \le n-1$, let $\iota_k^m : W_k^{(m)} \cong
W_k$ be an isomorphism to a copy of $W_k$.  Then for $1 \le k \le
n-1$, let
\begin{equation}
\label{eq:isom-map} \phi_k = \bigoplus_{p \in \mathcal{P},\,
\inc(p) = k} B_p i_{\out(p)} \iota_{\out(p)}^{\out(p)-\ord(p)} :
\bigoplus_{l=1}^{n-1} \bigoplus_{m\le k,l} W_l^{(m)} \to V_k.
\end{equation}

Let $d = \sum_{k=1}^{n-1} k w_k$ and identify $\bigoplus_{m,k:\, m
\le k} W_k^{(m)}$ with $\C^d$. Then $\theta : \sqcup_\v \L(\v,\w)
\to \F$ sends the point $[B,i,j]$ to the flag $F = (0 = F_0
\subset \dots \subset F_n = \C^d)$ where $F_k = \ker \phi_k$. Note
that $\theta$ is well-defined since the kernel of $\phi_k$ does
not change under the action of $G_\v$.  For $1 \le k \le n-1$, define
\[
W^{\le k} = \bigoplus_{m,l\, :\,  m \le l,k} W_l^{(m)}.
\]
Note that we always have
\[
F_k = \ker \phi_k \subset W^{\le k}.
\]

\begin{cor}
\label{cor:lang-isom} The image of the map $\theta : \sqcup_\v
\L(\v,\w) \to \F$ lies in $\F_x$ where $x \in N$ is the map
given in block form by
$W_k^{(m)} \stackrel{\cong}{\to} W_k^{(m-1)}$ (and
$x(W_k^{(1)})=0$). Furthermore $\theta : \sqcup_\v \L(\v,\w) \to
\F_x$ is an isomorphism and $\theta (\L(\v,\w)) = \F_{\a,x}$ where
$\a = \a(\v,\w)$ is defined by \eqref{eq:a-def}.
\end{cor}

\begin{prop}
\label{prop:ginz-res}
Let $\v,\w \in (\Z_{\ge 0})^{n-1}$, $\a=\a(\v,\w)$, $x \in N$ a nilpotent
element of type $1^{w_1} 2^{w_2} \cdots (n-1)^{w_{n-1}}$, and $1 \le k
\le n-1$.  Then
\begin{equation}
\label{eq:ginz-res} (\theta \times \theta) (\mathfrak{B}_k(\v,\w))
= \left( T^*_{Y_{\a_k^+,\a}}(\F_{\a_k^+} \times \F_\a) \right) \cap ({\tilde
S}_{\a_k^+,x} \times {\tilde S}_{\a,x}).
\end{equation}
\end{prop}

\begin{proof}
The right side of \eqref{eq:ginz-res} is equal to
\begin{equation} \label{eq:subflag}
\{(F',F) \in {\tilde S}_{\a_k^+,x} \times {\tilde S}_{\a,x}\ |\
F_l = F'_l\ \forall \ l \ne k,\, F_k \subset F'_k\,
\dim(F'_k/F_k) = 1 \}.
\end{equation}
Recall that
\[
\mathfrak{B}_k(\v,\w) = \{(B,i,j,S)\ |\ (B,i,j) \in
\mu^{-1}(0)^s,\, S \subset V,\, j(S)=0,\, S \text{ $B$-invariant},\, \dim
S=\e^k\}/G_\v.
\]
We consider this as a subset of $\M(\v-\e^k,\w) \times \M(\v,\w)$
by taking the quotient by the subspace $S$ in the first factor. We
know by Theorem~\ref{thm:maffei} that
\[
\theta : \M(\v,\w) \stackrel{\cong}{\longrightarrow} {\tilde
S}_{\a,x},\quad \theta : \M(\v-\e^k,\w)
\stackrel{\cong}{\longrightarrow} {\tilde S}_{\a_k^+,x}.
\]
Thus, it suffices to show that a choice of $B$-invariant subspace
$S$ of $V_k$ corresponds to a choice of $F'_k$ such that $F_k
\subset F'_k \subset x^{-1}(F_{k-1})$.  We first do this for the
case where $W = W_1$.  Then $i=i_1$ and $j=j_1$.
In this case, the isomorphism between
quiver varieties and flag varieties is particularly simple (see
\cite{N94} and \cite{M00}).  The isomorphism is given by $\theta :
[B,i,j] \mapsto (x,F)$ where
\begin{align*}
x=ji,\quad F=(0 \subset \ker i \subset \ker B_{12}i \subset \dots
\subset \ker B_{n-2,n-1} \dots B_{12} i \subset W).
\end{align*}
That is, $F_l = \ker B_{l-1,l} \cdots B_{12} i$.  Now, let $S
\subset V_k$ be a $B$-invariant subspace contained in the kernel
of $j$ with $\dim S = 1$ and let $(B',i',j')$ be the point of
$\M(\v-\e^k,\w)$ obtained from $(B,i,j)$ by taking the quotient by
the subspace $S$.  Now, since $S$ is $B$-invariant, we have that
$S \in \ker B_{k,k-1} \cap \ker B_{k,k+1}$.  Here we adopt the
convention that $B_{1,0}=0$ and $B_{n-1,n}=0$. Let $p : V_k \to
V_k/S$ be the canonical projection.  Then $\theta([B',i',j']) =
(x,F')$ where $x=ji$ and
\begin{align*}
F'_l &= \ker B_{l-1,l} \dots B_{12} i = F_l,\ l < k,\\
F'_l &= \ker B_{l-1,l} \dots B_{k,k+1} p B_{k-1,k} \dots B_{12}
i,\ l \ge k.
\end{align*}
Now, since $S \subset \ker B_{k,k+1}$, we have that $B_{k,k+1}p$ =
$B_{k,k+1}$.  Thus, for $l >k$, $F'_l = F_l$. Also,
\[
F_k' = \ker p B_{k-1,k} \cdots B_{12} i \supset \ker B_{k-1,k}
\cdots B_{12} i = F_k.
\]
Thus it remains to show
that $F'_k \subset x^{-1} (F_{k-1})$.  Now,
\begin{align*}
x^{-1}(F_{k-1}) &= x^{-1}(\ker B_{k-2,k-1} \dots B_{12} i) \\
&= \ker (B_{k-2,k-1} \dots B_{12} i x) \\
&= \ker (B_{k-2,k-1} \dots B_{12} iji).
\end{align*}
Now, since $(B,i,j) \in \mu^{-1}(0)$, we have that $ij =
B_{21}B_{12}$ and $B_{l-1,l} B_{l,l-1} = B_{l+1,l} B_{l+1,l}$ for
$2 \le l \le n-2$.  Thus,
\begin{align*}
B_{k-2,k-1} \dots B_{12} iji &= B_{k-2,k-1} \dots B_{12} B_{21}
B_{12} i \\
&\ \vdots \\
&=B_{k,k-1} B_{k-1,k} B_{k-2,k-1} \dots B_{12} i.
\end{align*}
Thus
\[
x^{-1}(F_{k-1}) = \ker (B_{k,k-1} B_{k-1,k} B_{k-2,k-1} \dots
B_{12} i).
\]
Now, since $S \subset \ker B_{k,k-1}$, we have
\[
F_k' = \ker(p B_{k-1,k} \dots B_{12}i) \subset \ker (B_{k,k-1}
B_{k-1,k} \dots B_{12} i) = x^{-1}(F_{k-1}).
\]
We have shown that every choice of subspace $S$ corresponds to a
flag $F'$ satisfying the conditions in \eqref{eq:subflag}.  It is
easy to see that such a flag $F'$ comes from a subspace $S$ as
follows.  We have that $F_k \subset F'_k$.  We take $S$ to
be the subspace of $V_k$ such that
\[
\ker (p B_{k-1,k} \dots B_{12} i) = F'_k
\]
for the projection $p : V_k \to V_k/S$.  Thus we have proven the
proposition in the special case $W = W_1$.

For the general case, we recall Maffei's construction in
\cite{M00}.  For general $W$, Maffei constructs a map
$\Lambda(\v,\w) \to \Lambda({\tilde \v},{\tilde \w})$, denoted
$(B,i,j) \mapsto ({\tilde B},{\tilde i},{\tilde j})$, where
${\tilde \w} = c\e^1$ for some $c \in \Z_{\ge 0}$. Thus, if we
show that a choice of a $B$-stable subspace $S$ such that $\dim S
= \e^k$ corresponds to a choice of $\tilde B$-stable subspace
$\tilde S$ such that $\dim {\tilde S} = \e^k$ then we reduce the
proof to the special case considered above.  Now,
\begin{align*}
{\tilde V}_k &= V_k \oplus W'_k, \\
\text{where } W'_k &= \bigoplus_{l,m\, :\, 1 \le m \le l-k,\, k+1
\le l \le n-1} W_l^{(m)},
\end{align*}
and $W_l^{(m)}$ is an isomorphic copy of $W_l$.  For $1 \le m \le
l-k$ and $k+1 \le l \le n-1$, we have (see \cite{M00})
\begin{gather*}
\pr_{W_l^{(m)}} {\tilde B}_{k,k-1} |_{W_l^{(m)}} = \Id_{W_l}, \\
\pr_{W_l^{(m)}} {\tilde B}_{k,k-1} |_{V_k} = 0,
\end{gather*}
where $\pr_{W_l^{(m)}}$ denotes the projection onto the subspace
$W_l^{(m)}$.   In particular, $\ker {\tilde B}_{k,k-1} \subset
V_k$. Thus, since the subspace ${\tilde S} \subset {\tilde V}_k$
must be contained in $\ker {\tilde B}_{k,k-1}$, it must lie in
$V_k$. The result then follows from Remark~19 of \cite{M00}.
\end{proof}

We now compare the Lie algebra action in the two settings.  By
\cite[\S 3.7.14]{CG}, $S_x$ is transverse to the orbit $O_x$ in
$N$. Thus, there is an open neighborhood $U \subset N$ of $S$
such that
\[
U \cong (O_x \cap U) \times S.
\]
Let $\tilde U_\d = \mu_\d^{-1}(U)$ and $M_\d' = \mu_\d^{-1}(S_x) =
\mu_\d^{-1}(S_{\d,x}) = \tilde S_{\d,x}$.
Then $\tilde U_\d \subset M_\d$ is an open
neighborhood of $M_\d'$. Let $D = O_x \cap U$, a small
neighborhood of $x$ in $O_x$.  By \cite[Cor 3.2.21]{CG},
\[
\tilde U_\d \cong (O_x \cap U) \times M_\d'.
\]
Then the two commutative diagrams
\[
\begin{CD}
M_\d' @>>> M_\d \\
@VV{\mu_\d}V  @VV{\mu_\d}V \\
S_x @>>> N
\end{CD}
\qquad \text{and} \qquad
\begin{CD}
M_\d' @>{p \mapsto (x,p)}>> D \times M_\d' @>{\cong}>> \tilde U_\d \\
@VV{\mu_\d}V  @VV{\id_D \times \mu_\d}V  @VV{\mu_\d}V \\
S_x @>{y \mapsto (x,y)}>> D \times S_x @>{\cong}>> U
\end{CD}
\]
are isomorphic, where the horizontal arrows in the left diagram
are given by the natural inclusions.  If we let $\tilde U =
\mu^{-1}(U)$ and $M' = \mu^{-1}(S_x)$ then $\tilde U = \sqcup_\d
\tilde U_\d$, $M' = \sqcup_\d M_\d'$ and $\tilde U \subset M$ is
an open neighborhood of $M'$.  Thus we have that the two
commutative diagrams
\[
\begin{CD}
M' @>>> M \\
@VV{\mu}V  @VV{\mu}V \\
S_x @>>> N
\end{CD}
\qquad \text{and} \qquad
\begin{CD}
M' @>{p \mapsto (x,p)}>> D \times M' @>{\cong}>> \tilde U \\
@VV{\mu}V  @VV{\id_D \times \mu}V  @VV{\mu}V \\
S_x @>{y \mapsto (x,y)}>> D \times S_x @>{\cong}>> U
\end{CD}
\]
are isomorphic.  Let $Z' = M' \times_{S_x} M'$.  Then by
Theorem~\ref{thm:maffei},
\[
Z' \cong \sqcup_{\v^1,\v^2} Z(\v^1,\v^2;\w).
\]
We then have the commutative diagram
\begin{equation}
\label{eq:comm-slice}
\begin{CD}
Z' = M' \times_{S_x} M' @>i>> M \times_N M = Z \\
@VVV  @VVV \\
M' \times M' @>i>> M \times M
\end{CD}
\end{equation}
where the maps are the obvious inclusions.
Diagram~\eqref{eq:comm-slice} is isomorphic to
\begin{equation}
\begin{CD}
Z' = M' \times_{S_x} M' @>i={p \mapsto (x,p)}>> D_\Delta \times
(M' \times_{S_x} M')
\cong Z \cap (\tilde U \times \tilde U) \\
@V{j}VV  @V{\Delta \times j}VV \\
M' \times M' @>i={p \mapsto ((x,x),p)}>> (D \times D) \times (M'
\times M')
\end{CD}
\end{equation}
where $\Delta : D_\Delta \to D \times D$ is the embedding of the
diagonal.

\begin{lem}
\label{lem:irrcomp-restrict}
The inverse image in $M' \times M'$ of an irreducible component of
the variety $Z$ is either empty or else is an irreducible
component of the variety $Z'$.
\end{lem}
\begin{proof}
Let $X$ be a (closed) irreducible component of $Z$.  If $X$ does
not intersect the open subset $\tilde U \times \tilde U \subset M
\times M$, then $i^{-1}(X) = \emptyset$, since $i(Z') \subset
\tilde U \times \tilde U$.  Now assume that $X_U = X \cap (\tilde
U \times \tilde U)$ is non-empty.  Then $X_U$ is an irreducible
component of $Z \cap (\tilde U \times \tilde U)$. Thus is must be
of the form $X_\U \cong D_\Delta \times X'$ where $X'$ is an
irreducible component of $M' \times_{S_x} M' = Z'$.  We then have
$i^{-1}(X) = X'$ and the result follows.
\end{proof}

The diagram~\eqref{eq:comm-slice} gives rise to a restriction with
support morphism
\[
i^* : H_*(Z) \to H_*(Z'),\quad c \mapsto c \cap [M' \times M'].
\]
By Lemma~\ref{lem:irrcomp-restrict}, $i^*$ takes $H_\t(Z)$
to $H_\t(Z')$.  Furthermore, by
Proposition~\ref{prop:ginz-res} we have that
\begin{equation}
\label{eq:action-isom}
i^*([T^*_{Y_{\a_k^+,\a}}(\F_{\a_k^+} \times \F_\a)]) = [(\theta
\times \theta)(\mathfrak{B}_k(\v,\w))]
\end{equation}
where $\a = \a(\v,\w)$.

Now, $\F_x = \mu^{-1}(x)$ can be viewed
as a subvariety of $M'$ or $M$.  If $i : M' \to M$ is the
inclusion, then the restriction with supports morphism $i^* :
H_\t(\F_x) \to H_\t(\F_x)$ is an isomorphism, where the first and second
$H_\t(\F_x)$ are $H^0(M,M\backslash {\F_x})$ and $H^0(M',M'\backslash
{\F_x})$ respectively.

\begin{theo}
\label{thm:quiver-flag}
\begin{enumerate}
\item The morphism $i^* : H_\t(Z) \to H_\t(Z') \cong \oplus_{\v^1,\v^2}
H_\t(Z(\v^1,\v^2;\w))$ is an algebra homomorphism (with respect to
the convolution product).

\item The following diagram, where $x \in N$ is a nilpotent element of type
$1^{w_1} 2^{w_2} \cdots (n-1)^{w_{n-1}}$ and whose vertical maps are given by
convolution, commutes
\[
\begin{CD}
H_\t(Z) \otimes H_\t(\F_x) @>{i^* \otimes i^*}>> H_\t(Z') \otimes
H_\t(\F_x) @>\cong>> \oplus_{\v^1,\v^2} H_\t(Z(\v^1,\v_2;\w))
\otimes \oplus_\v H_\t(\L(\v,\w))
\\
@VVV  @VVV  @VVV \\
H_\t(\F_x) @>{i^*}>> H_\t(\F_x) @>\cong>> \oplus_\v
H_\t(\L(\v,\w))
\end{CD}
\]
\end{enumerate}
\end{theo}
\begin{proof}
Note that the two rightmost horizontal maps are the isomorphisms induced
by the map $\theta$ of Theorem~\ref{thm:maffei}.
We prove only the first part of the theorem.  The second part in
analogous.  We have a sequence of embeddings
\[
M' \times M' \hookrightarrow \tilde U \times \tilde U
\hookrightarrow M \times M
\]
So $i^*$ factors as
\[
i^* : H_\t(Z) \to H_\t(Z \cap (\tilde U \times \tilde U)) \to H_\t(Z')
\]
The first map is the restriction to an open subset and thus
commutes with convolution by base locality (cf. \cite[\S
2.7.45]{CG}).  The second map is induced by the embedding
\[
Z' \hookrightarrow Z \cap (\tilde U \times \tilde U).
\]
By the above results, this is isomorphic to the natural embedding
\[
Z' \hookrightarrow D_\Delta \times Z',\quad z \mapsto (x,z).
\]
The corresponding map
\[
i^* : H_\t(D_\Delta \times Z') \to H_\t(Z')
\]
commutes with convolution by the K\" unneth formula for
convolution (cf. \cite[\S 2.7.16]{CG}).
\end{proof}

\begin{cor}
\label{cor:ginz=nak}
If $c \in H_\t(F_x)$ and $c'$ is the corresponding class in
$\oplus_\v H_\t(\L(\v,\w))$ (under the isomorphism $\theta$) then
we have
\begin{align*}
E^\text{Gin}_k c &= E^\text{Nak}_k c', \text{ and} \\
F^\text{Gin}_k c &= F^\text{Nak}_k c' \text{ for all $k$}.
\end{align*}
Here the superscripts Gin and Nak correspond to the actions
defined by Ginzburg and Nakajima respectively.
\end{cor}
\begin{proof}
The result follows from \eqref{eq:action-isom} and the fact that
since $\tilde U_\d \cong (O_x \cap U) \times M_\d'$ we have
\begin{align*}
\dim_\C M_{\d^1} - \dim_\C M_{\d^2} &= (\dim_\C (O_x \cap U) +
\dim_\C M_{\d^1}') -( \dim_\C (O_x \cap U) + \dim_\C M_{\d^2}') \\
&= \dim_\C M_{\d^1}' - \dim_\C M_{\d^2}'.
\end{align*}
Thus the signs appearing in \eqref{def:ginz-Fk} and
\eqref{def:nak-Fk} are the same.
\end{proof}

We see from Corollary~\ref{cor:ginz=nak} that the Ginzburg and
Nakajima constructions yield the same representations, with the
same bases, given by the fundamental classes of the irreducibles
components of $\F_x \cong \sqcup_\v \L(\v,\w)$. However, note that
the corresponding quotients of the universal enveloping algebra
constructed via convolution is different (compare
Theorems~\ref{thm:ginz-quotient} and \ref{thm:nak-quotient}).  To
see that these two quotients are indeed different, it suffices to
consider the case of $\mathfrak{sl}_3$ with $\w = (1,1)$ (so
$\omega_\w = \omega_1 + \omega_2$ and $d=3$).  Then the weight
$3\omega_1$ corresponds to a partition of $d$ but is not a weight
of $L(\omega_\w)$ (since the tableau of shape $(21)$ with all
three entries equal to 1 is not semistandard).

%%%%%%%%%%%%%%%%%%%%%%%%%%%%%%%%%%%%%%%%%%%%%%%%%%%%%%%%%%%%%%%%%%%%%%%

\section{Crystal structure on flag varieties}
\label{sec:crystal}

Kashiwara and Saito have introduced the structure of a crystal on
the set of irreducible components of Nakajima's quiver varieties. In
this section, we recall this construction and use the isomorphism of
Section~\ref{sec:isom} to define a crystal structure on the flag
varieties (or, more precisely, on the set of irreducible components
of the Spaltenstein varieties $\F_x$). In this way we recover the
crystal structure defined by Malkin (see \cite{Mal02}). In fact,
Malkin and Nakajima have defined a tensor product quiver variety
(see \cite{Mal03} and \cite{Nak01}). One would expect that the
relationship between the two constructions examined in this paper
could be extended to this setting and one would recover the tensor
product crystal structure defined in \cite{Mal02}.  However, we will
restrict ourselves to the case of a single representation here.

We first review the realization of the crystal graph via quiver
varieties. See \cite{KS97,S02} for proofs omitted here. Note that,
as mentioned in Section~\ref{sec:nak}, we are using a different
stability condition and thus our definitions differ slightly from
those in \cite{KS97,S02}.

Let $\mathbf{w, v, v', v''} \in (\Z_{\ge 0})^I$ be such that $\v =
\mathbf{v'} + \mathbf{v''}$.  Consider the maps
\begin{equation}
\label{eq:diag_action} \Lambda(\v'',\mathbf{0}) \times
\Lambda(\v',\w) \stackrel{p_1}{\longleftarrow} \mathbf{\tilde F
(v,w;v'')} \stackrel{p_2}{\longrightarrow} \mathbf{F(v,w;v'')}
\stackrel{p_3}{\rightarrow} \Lambda(\v,\w),
\end{equation}
where the notation is as follows.  A point of
$\mathbf{F(v,w;v'')}$ is a point $(B,i) \in \Lambda(\v,\w)$
together with an $I$-graded, $B$-stable subspace $S$ of $V$ such
that $\dim S = \v''$.  A point of $\mathbf{\tilde
  F (v,w;v'')}$ is a point $(B,i,S)$ of $\mathbf{F(v,w;v'')}$
together with a collection of isomorphisms $R''_k : V''_k \cong
S_k$ and $R'_k : V'_k \cong V_k / S_k$ for each $k \in I$. Then we
define $p_2(B,i,S, R',R'') = (B,i,S)$, $p_3(B,i,S) = (B,i)$ and
$p_1(B,i,S,R',R'') = (B'',B',i')$ where $B'', B', i'$ are
determined by
\begin{align*}
R''_{\inc(h)} B''_h &= B_h R''_{\out(h)} : V''_{\out(h)} \to
S_{\inc(h)}, \\
R'_k i'_k &= {\bar i}_k : W_k \to V_k/S_k, \\
R'_{\inc(h)} B'_h &= B_h R'_{\out(h)} : V'_{\out(h)} \to
V_{\inc(h)}/S_{\inc(h)},
\end{align*}
where ${\bar i}_k$ denotes the composition of the map $i_k$ with
the canonical projection $V_k \to V_k/S_k$. It follows that $B'$
and $B''$ are nilpotent.

\begin{lem}[{\cite[Lemma 10.3]{N94}}]
One has
\[
(p_3 \circ p_2)^{-1} (\Lambda(\v,\w)^s) \subset p_1^{-1}
(\Lambda(\v'',\mathbf{0}) \times \Lambda(\v',\w)^s).
\]
\end{lem}

Thus, we can restrict \eqref{eq:diag_action} to stable points,
forget the $\Lambda(\v'',\mathbf{0})$-factor and consider the
quotient by $G_\v$, $G_{\v'}$. This yields the diagram
\begin{equation}
\label{eq:diag_action_mod} \mathcal{L}(\v', \w)
\stackrel{\pi_1}{\longleftarrow} \mathcal{L}(\v, \w; \v -
\mathbf{v'}) \stackrel{\pi_2}{\longrightarrow} \mathcal{L}(\v,
\w),
\end{equation}
where
\[
\mathcal{L}(\v, \w; \v - \mathbf{v'}) \stackrel{\text{def}}{=} \{
(B,i,S) \in \mathbf{F(\v,\w;\v-\v')}\,
  |\, (B,i) \in \Lambda(\v,\w)^s \} / G_\v.
\]

For $k \in I$ define $\varepsilon_k : \Lambda(\v, \w) \to \Z_{\ge
0}$ by
\[
\varepsilon_k((B,i)) = \dim_\C \ker \left( V_k
\stackrel{(B_h)}{\longrightarrow} \bigoplus_{h\, :\, \out(h)=k}
V_{\inc(h)} \right).
\]
Then, for $c \in \Z_{\ge 0}$, let
\[
\mathcal{L}(\v,\w)_{k,c} = \{[B,i] \in \mathcal{L}(\v,\w)\ |\
\varepsilon_k((B,i)) = c\}
\]
where $[B,i]$ denotes the $G_\v$-orbit through the point $(B,i)$.
$\mathcal{L}(\v,\w)_{k,c}$ is a locally closed subvariety of
$\mathcal{L}(\v,\w)$.

Assume $\mathcal{L}(\v,\w)_{k,c} \ne \emptyset$ and let $\v' = \v
- c\mathbf{e}^k$ where $\mathbf{e}^k_l = \delta_{kl}$.  Then
\[
\pi_1^{-1}(\mathcal{L}(\v',\w)_{k,0}) =
\pi_2^{-1}(\mathcal{L}(\v,\w)_{k,c}).
\]
Let
\[
\mathcal{L}(\v,\w;c\mathbf{e}^k)_{k,0} =
\pi_1^{-1}(\mathcal{L}(\v',\w)_{k,0}) =
\pi_2^{-1}(\mathcal{L}(\v,\w)_{k,c}).
\]
We then have the following diagram.
\begin{equation}
\label{eq:crystal-action} \mathcal{L}(\v',\w)_{k,0}
\stackrel{\pi_1}{\longleftarrow} \mathcal{L}(\v,\w;
c\mathbf{e}^k)_{k,0} \stackrel{\pi_2}{\longrightarrow}
\mathcal{L}(\v,\w)_{k,c}
\end{equation}
The restriction of $\pi_2$ to $\mathcal{L}(\v,\w;
c\mathbf{e}^k)_{k,0}$ is an isomorphism since the only possible
choice for the subspace $S$ of $V$ is to have $S_l = 0$ for $l \ne
k$ and $S_k$ equal to the intersection of the kernels of the $B_h$
with $\out(h)=k$. Also, $\mathcal{L}(\v',\w)_{k,0}$ is an open
subvariety of $\mathcal{L}(\v',\w)$.

\begin{lem}[\cite{S02}]
\begin{enumerate}
\item For any $k \in I$,
\[
\mathcal{L}(\mathbf{0},\w)_{k,c} =
\begin{cases}
pt & \text{if $c=0$} \\
\emptyset & \text{if $c >0$}
\end{cases}.
\]
\item Suppose $\mathcal{L}(\v,\w)_{k,c} \ne \emptyset$ and $\v' =
\v - c\mathbf{e}^k$.  Then the fiber of the restriction of $\pi_1$
to $\mathcal{L}(\v, \w; c\mathbf{e}^k)_{k,0}$ is isomorphic to a
Grassmannian variety.
\end{enumerate}
\end{lem}

\begin{cor}
\label{cor:irrcomp-isom} Suppose $\mathcal{L}(\v,\w)_{k,c} \ne
\emptyset$.  Then there is a 1-1 correspondence between the set of
irreducible components of $\mathcal{L}(\v - c\mathbf{e}^k,
\w)_{k,0}$ and the set of irreducible components of
$\mathcal{L}(\v, \w)_{k,c}$.
\end{cor}

Let $\mathcal{B}(\v,\w)$ denote the set of irreducible components
of $\mathcal{L}(\v,\w)$ and let $\mathcal{B}(\w) = \bigsqcup_\v
\mathcal{B}(\v,\w)$. For $X \in \mathcal{B}(\v,\w)$, let
$\varepsilon_k(X) = \varepsilon_k((B,i))$ for a generic point
$[B,i] \in X$.  Then for $c \in \Z_{\ge 0}$ define
\[
\mathcal{B}(\v,\w)_{k,c} = \{X \in \mathcal{B}(\v,\w)\ |\
\varepsilon_k(X) = c\}.
\]
Then by Corollary~\ref{cor:irrcomp-isom}, $\mathcal{B}(\v -
c\mathbf{e}^k,\w)_{k,0} \cong \mathcal{B}(\v, \w)_{k,c}$.

Suppose that ${\bar X} \in \mathcal{B}(\v -
c\mathbf{e}^k,\w)_{k,0}$ corresponds to $X \in
\mathcal{B}(\v,\w)_{k,c}$ by the above isomorphism. Then we define
maps
\begin{gather*}
\kf_k^c : \mathcal{B}(\v - c\mathbf{e}^k,\w)_{k,0} \to
\mathcal{B}(\v,\w)_{k,c},\quad \kf_k^c({\bar X})
= X, \\
\ke_k^c : \mathcal{B}(\v,\w)_{k,c} \to \mathcal{B}(\v -
c\mathbf{e}^k,\w)_{k,0},\quad \ke_k^c(X) = {\bar X}.
\end{gather*}
We then define the maps
\[
\ke_k, \kf_k : \mathcal{B}(\w) \to \mathcal{B}(\w) \sqcup \{0\}
\]
by
\begin{gather*}
\ke_k : \mathcal{B}(\v,\w)_{k,c}
\stackrel{\ke_k^c}{\longrightarrow} \mathcal{B}(\v -
c\mathbf{e}^k, \w)_{k,0} \stackrel{\kf_k^{c-1}}{\longrightarrow}
\mathcal{B}(\v -
\mathbf{e}^k, \w)_{k,c-1}, \\
\kf_k : \mathcal{B}(\v,\w)_{k,c}
\stackrel{\ke_k^c}{\longrightarrow} \mathcal{B}(\v -
c\mathbf{e}^k, \w)_{k,0} \stackrel{\kf_k^{c+1}}{\longrightarrow}
\mathcal{B}(\v + \mathbf{e}^k, \w)_{k,c+1}.
\end{gather*}
We set $\ke_k(X)=0$ for $X \in \mathcal{B}(\v,\w)_{k,0}$ and
$\kf_k(X)=0$ for $X \in \mathcal{B}(\v,\w)_{k,c}$ with
$\mathcal{B}(\v,\w)_{k,c+1} = \emptyset$. We also define
\begin{gather*}
\wt : \mathcal{B}(\w) \to P,\quad \wt(X) = \omega_\w - \alpha_\v
\text{
  for } X \in \mathcal{B}(\v,\w), \\
\varphi_k(X) = \varepsilon_k(X) + \left< h_k, \wt(X) \right>.
\end{gather*}

\begin{prop}[\cite{S02}]
\label{prop:quiver-crystal}
$\mathcal{B}(\w)$ is a crystal and is isomorphic to the crystal of
the highest weight $U_q(\g)$-module with highest weight
$\omega_\w$.
\end{prop}

We now translate this structure to the language of flag varieties.
We need the following results.  We adopt the convention that
$B_{1,0}=0$ and $B_{n-1,n}=0$.

\begin{prop}
\label{prop:comm1} We have
\[
B_{k,k-1} \circ \phi_k = \phi_{k-1} \circ x.
\]
\end{prop}

\begin{proof}
Recall that
\begin{align*}
\phi_{k-1} &= \bigoplus_{p \in \mathcal{P},\,
\inc(p) = k-1} B_p i_{\out(p)} \iota_{\out(p)}^{\out(p)-\ord(p)} \\
\Rightarrow \phi_{k-1} \circ x &= \bigoplus_{p \in \mathcal{P},\,
\inc(p) = k-1,\, \ord(p) \ge 1} B_p i_{\out(p)}
\iota_{\out(p)}^{\out(p)-\ord(p)+1}.
\end{align*}
And
\begin{align*}
\phi_k &= \bigoplus_{p \in \mathcal{P},\,
\inc(p) = k} B_p i_{\out(p)} \iota_{\out(p)}^{\out(p)-\ord(p)} \\
\Rightarrow B_{k,k-1} \circ \phi_k &= \bigoplus_{p \in
\mathcal{P},\, \inc(p) = k} B_{k,k-1} B_p i_{\out(p)}
\iota_{\out(p)}^{\out(p)-\ord(p)}.
\end{align*}
Now, since $j=0$, $\mu(B,i,j)=0$ implies that
\begin{gather*}
B_{l-1,l} B_{l,l-1} = B_{l+1,l} B_{l,l+1} \text{ for } 2 \le l \le
n-2, \\
B_{2,1} B_{1,2} = 0, \quad B_{n-2,n-1} B_{n-1,n-2} = 0.
\end{gather*}
Using these equations, one can see that
\[
\{B_{k,k-1} B_p\ |\  p \in \P,\, \inc(p)=k\} = \{B_p\ |\ p \in
\P,\, \inc(p) = k-1,\, \ord(p) \ge 1\}.
\]
Therefore,
\[
B_{k,k-1} \circ \phi_k = \bigoplus_{p \in \P,\, \inc(p)=k-1,\,
\ord(p) \ge 1} B_p i_{\out(p)}
\iota_{\out(p)}^{\out(p)-(\ord(p)-1)} = \phi_{k-1} \circ x.
\]
\end{proof}

\begin{prop}
\label{prop:comm2} We have
\[
B_{k,k+1} \circ \phi_k = \phi_{k+1} |_{W^{\le k}}.
\]
\end{prop}

\begin{proof}
Let $\P'$ be the subset of $\P$ consisting of those paths that
contain at least one edge belonging to $\bar \Omega$.  Then
\begin{align*}
B_{k,k+1} \circ \phi_k &= \bigoplus_{p \in \mathcal{P},\,
\inc(p) = k} B_{k,k+1} B_p i_{\out(p)} \iota_{\out(p)}^{\out(p)-\ord(p)} \\
&= \bigoplus_{p \in \mathcal{P}',\,
\inc(p) = k+1} B_p i_{\out(p)} \iota_{\out(p)}^{\out(p)-\ord(p)} \\
&= \phi_{k+1} |_{W^{\le k}}.
\end{align*}
\end{proof}

\begin{prop}
\label{prop:flag-subspace}
One has
\begin{equation}
\label{eq:flag-subspace} \phi_k^{-1} \left( \ker B_{k,k-1} \cap
\ker B_{k,k+1} \right) = x^{-1}(F_{k-1}) \cap F_{k+1}.
\end{equation}
\end{prop}

\begin{proof}
Since $F_{k-1} \subset W^{\le k-1}$, we have that $x^{-1}(F_{k-1})
\subset W^{\le k}$.  Thus, using propositions~\ref{prop:comm1} and
\ref{prop:comm2},
\begin{align*}
x^{-1}(F_{k-1}) \cap F_{k+1} &= x^{-1}(F_{k-1}) \cap
\ker(\phi_{k+1}) \\
&= x^{-1}(F_{k-1}) \cap \ker(\phi_{k+1}|_{W^{\le k}}) \\
&= x^{-1}(\ker \phi_{k-1}) \cap \ker(\phi_{k+1}|_{W^{\le k}}) \\
&= \ker (\phi_{k-1}\circ x) \cap \ker(\phi_{k+1}|_{W^{\le k}}) \\
&= \ker (B_{k,k-1} \circ \phi_k) \cap \ker (B_{k,k+1} \circ
\phi_k) \\
&= \phi_k^{-1} \left( \ker B_{k,k-1} \cap \ker B_{k,k+1} \right)
\end{align*}
\end{proof}

Note that
\begin{equation}
\label{eq:stable-subspace} \ker B_{k,k-1} \cap \ker B_{k,k+1} =
\ker \left( V_k \stackrel{(B_h)}{\longrightarrow} \bigoplus_{h\,
:\, \out(h)=k} V_{\inc(h)} \right),
\end{equation}
and that a collection of subspaces $S_l \subset V_l$ such that
$S_l = 0$ for $l \ne k$ is $B$-stable if and only if $S_k$ is
contained in the righthand side of
equation~\eqref{eq:stable-subspace}. Thus, the flag variety
analogue of the diagram~\eqref{eq:crystal-action} (for $\v - \v' =
c\mathbf{e}^k$) is
\[
\F_{\a^{k,c},x} \stackrel{\pi_1}{\longleftarrow} \F_{\a,x}(k,c)
\stackrel{\pi_2}{\longrightarrow} \F_{\a,x},
\]
where $\a = \a(\v,\w)= (a_1,\dots,a_n)$ and $\a^{k,c} =
(a_1,\dots,a_{k-1},a_k+c, a_{k+1}-c,a_{k+2},\dots,a_n)$, and
\[
\F_{\a,x}(k,c) = \{(F,S)\ |\ F \in \F_{\a,x},\, F_k \subset S
\subset F_{k+1} \cap x^{-1}(F_{k-1}),\, \dim S/F_k = c\}.
\]
In particular
\[
\mathcal{L}(\v,\w;c\mathbf{e}^k) \cong \F_{\a,x}(k,c).
\]

Let $\mathcal{B}(\a,x)$ denote the set of irreducible components
of $\F_{\a,x}$ and let $\mathcal{B}(x) = \sqcup_\a \F_{\a,x}$. Let
\[
\varepsilon_k(F) = \dim (F_{k+1} \cap x^{-1}(F_{k-1})) - \dim F_k,
\]
and for $X \in \mathcal{B}(\a,x)$ define $\varepsilon_k(X) =
\varepsilon_k(F)$ for a generic flag $F \in X$.  Then for $c \in
\Z_{\ge 0}$ define
\[
\mathcal{B}(\a,x)_{k,c} = \{X \in \mathcal{B}(\a,x)\ |\
\varepsilon_k(X)=c\}.
\]
Then just as for quiver varieties, we have
$\mathcal{B}(\a^{k,c},x)_{k,0} \cong \mathcal{B}(\a,x)_{k,c}$ and
we define $\kf_k$ and $\ke_k$ just as before. We also define
\begin{gather*}
\wt(X) : \mathcal{B}(x) \to P,\quad \wt(X) = \sum_{k \in I} a_k
\epsilon_k
\text{ for } X \in \mathcal{B}(\a,x), \\
\varphi_k(X) = \varepsilon_k(X) + \left<h_k,\wt(X)\right>.
\end{gather*}
Then, by translating Proposition~\ref{prop:quiver-crystal} into
the language of flag varieties, we have the following theorem.

\begin{theo}
$\mathcal{B}(x)$ is a crystal and is isomorphic to the crystal of
the highest weight $U_q(\mathfrak{sl}_n)$-module with highest
weight $w_1 \omega_1 + \dots + w_{n-1} \omega_{n-1}$ where
$\omega_i$ are the fundamental weights of $\mathfrak{sl}_n$ and
$w_i$ is the number of $(i \times i)$-Jordan blocks in the Jordan
normal form of $x$.
\end{theo}

%%%%%%%%%%%%%%%%%%%%%%%%%%%%%%%%%%%%%%%%%%%%%%%%%%%%%%%%%%%%%%%%%%%%%%%

\bibliographystyle{abbrv}
\bibliography{biblist}

\end{document}